\documentclass[letterpaper, 10 pt, conference]{ieeeconf}
\usepackage{generic}
\usepackage{cite}
\usepackage{amsmath,amssymb,amsfonts}
\usepackage{graphicx}
\usepackage{textcomp}
\usepackage{cite}

\usepackage{bbold}
\usepackage{amsmath,amssymb,amsfonts}
\usepackage{graphicx}
\usepackage{textcomp}
\usepackage{lipsum}

\usepackage{textcomp}
\usepackage{subcaption}

\usepackage{verbatim}

\usepackage{algorithm}
\usepackage{algpseudocode}

\usepackage{hyperref}

\usepackage{amssymb}
\usepackage{amsmath}
\usepackage{multirow}
\usepackage{adjustbox}
\usepackage{booktabs}
\usepackage{tikz}
\usepackage{pgfplots}
\pgfplotsset{compat=1.16} 
\usepackage{tikz}
\usetikzlibrary{calc}
\usetikzlibrary{decorations.pathmorphing}
\usetikzlibrary{arrows, arrows.meta}
\usetikzlibrary{calc,shapes.geometric,arrows,positioning,intersections}
\usetikzlibrary{decorations, decorations.text,backgrounds}
\newtheorem{lemma}{Lemma}
\newtheorem{corollary}{Corollary}
\newtheorem{definition}{Definition}
\newtheorem{assumption}{Assumption}
\newtheorem{theorem}{Theorem}

\DeclareMathOperator{\st}{s.t.:}

\usepackage{titlesec}
\titlespacing*{\section}
{0pt} 
{1pt} 
{1pt} 

\titlespacing*{\subsection}
{0pt} 
{0.5pt} 
{0.5pt} 

\newcommand{\cP}{{\cal P}}
\newcommand{\cN}{{\cal N}}

\def\BibTeX{{\rm B\kern-.05em{\sc i\kern-.025em b}\kern-.08em
    T\kern-.1667em\lower.7ex\hbox{E}\kern-.125emX}}
\begin{document}
\title{A-Priori Reduction of Scenario Approximation for Automated Generation Control in High-Voltage Power Grids with Renewable Energy
}
\author{Aleksandr Lukashevich, 
Aleksandr Bulkin, and Yury Maximov, \IEEEmembership{Senior Member, IEEE}
\thanks{A. Lukashevich and A. Bulkin are with the Skolkovo Institue of Science and Technology,  
Moscow, 125043, and Institute for Artificial Intelligence, Lomonosov Moscow State University, Moscow, 119191, Russia. A. Bulkin is also with the International Center of Corporate Data Analysis, Astana, 010010, Kazakhstan. Y. Maximov is with the Los Alamos National Laboratory, 87545, Los Alamos, NM, USA.
}
\vspace{-0mm}
}
\maketitle
\thispagestyle{empty}

\begin{abstract}

Renewable energy sources (RES) are increasingly integrated into power systems to support the United Nations' Sustainable Development Goals of decarbonization and energy security. However, their low inertia and high uncertainty pose challenges to grid stability and increase the risk of blackouts. Stochastic chance-constrained optimization, particularly data-driven methods, offers solutions but can be time-consuming, especially when handling multiple system snapshots. This paper addresses a dynamic joint chance-constrained Direct Current Optimal Power Flow (DC-OPF) problem with Automated Generation Control (AGC) to facilitate cost-effective power generation while ensuring that balance and security constraints are met. We propose an approach for a data-driven approximation that includes a priori sample reduction, maintaining solution reliability while reducing the size of the data-driven approximation. Both theoretical analysis and empirical results demonstrate the superiority of this approach in handling generation uncertainty, requiring up to twice less data while preserving solution reliability.
\end{abstract}

\vspace{-0mm}
\section{Introduction}
\label{sec:introduction}
\vspace{-0mm}
Integrating renewable energy sources align with the United Nations' sustainable development goals, promoting affordable and clean energy while enhancing energy security and resilience. Unfortunately, RES introduces significant uncertainty in power systems generation, posing substantial challenges to grid optimization and control policies.

The Optimal Power Flow (OPF) \cite{stott2012optimal} is a key optimization problem that aims to achieve economically optimal generation while adhering to grid security and power balance constraints. To account for generation uncertainty, the Joint Chance-Constrained (JCC) extension considers an unknown joint distribution of renewable energy sources \cite{geng2019data, bienstock2014chance}. An alternative robust optimization approach assumes bounded uncertainty and offers a more conservative solution in practice \cite{ben2002robust, ding2016adjustable}.
The discrete-time dynamic chance-constrained OPF problem \cite{lou2019multi, monticelli1987security} models optimal generation set-points for sequential timestamps, temporarily binding generators' power outputs through ramp-up and ramp-down constraints. These constraints model the limit of the rate of change of the power output, as significant immediate changes are not feasible for some generators \cite{frangioni2008solving}. Automatic Generation Control (AGC) is widely used for fast and efficient power dispatch in bulk power systems \cite{xu2017real}.


While the chance-constrained extension enhances flexibility in modeling uncertainty, solving it for an arbitrary distribution and/or jointly for all technical limits becomes computationally infeasible~\cite{nemirovski2012safe}. To overcome this, data-driven (DD) approximations such as Scenario Approximation (SA) \cite{calafiore2006scenario} and Sample Average Approximation (SAA)~\cite{ahmed2008solving} have proven successful.
{\color{black}{}
Data-driven scenario approaches are often computationally prohibitive when higher accuracy is needed, prompting extensive scenario reduction studies. Scenario reduction methods are either a-posteriori or a-priori. A-posteriori methods solve the initial JCC problem multiple times to identify reducible scenarios iteratively \cite{campi2011sampling, geng2019data}. A-priori methods reduce scenarios before solving, enhancing computational efficiency. These methods, pioneered in~\cite{dupavcova2003scenario, dupavcova1990stability} and refined in~\cite{heitsch2003scenario}, use probability metrics like the Wasserstein distance to drop scenarios. Alternatively, some methods use clustering to replace the initial set with a smaller, more representative set \cite{rujeerapaiboon2022scenario, keutchayan2023problem}.}
{\color{black}{}The key difference between a-posteriori and a-priori reductions is the theoretical guarantees for SA solution feasibility for the original JCC problem. To our knowledge, such guarantees are limited for a-priori methods in the literature. To address this, we propose the A-priori Reduced Scenario Approximation (AR-SA), an approach for a-priori sample reduction methods linked with a data-driven Scenario Approximation (SA). This SA requires significantly fewer samples to produce a reliable solution for JCC dynamic DC optimal power flow and provides theoretical feasibility guarantees for JCC DC-OPF.}

\begin{figure}
    \centering
    \hspace{-2mm}\includegraphics[width=0.43\textwidth]{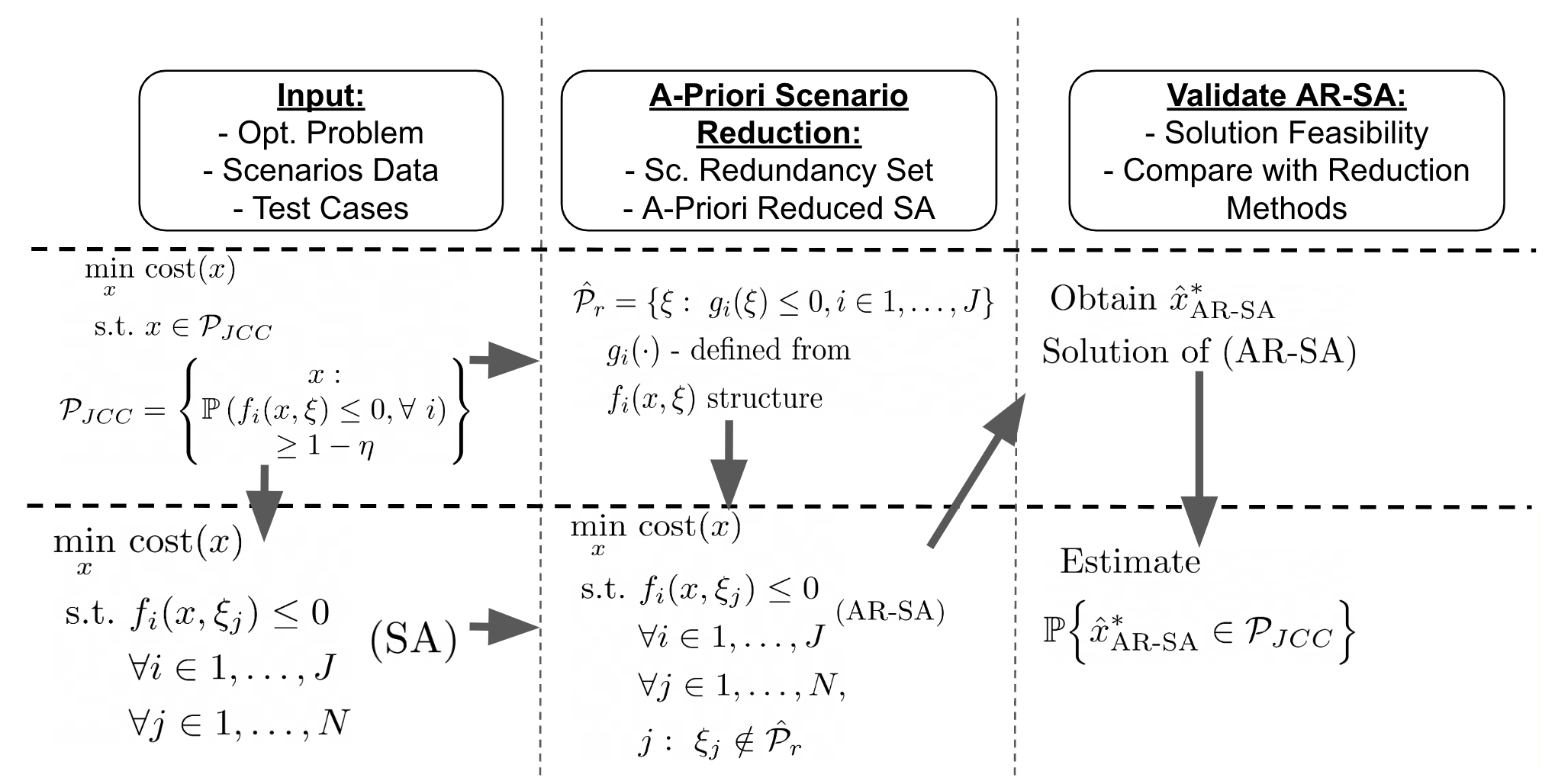}
    \caption{AR-SA workflow. We compare the performance of AR-SA with other reduction techniques such as Fast Forward, Simultaneous Backward, and K-Means methods on Grid6-WW, Washington-14, and IEEE-30 grids.}
    \label{fig:workflow}
    \vspace{-0mm}
\end{figure}

The contributions of this paper are as follows. First, we analytically define a-priori conditions that determine sample redundancy for JCC dynamic optimal power flow with AGC and provide theoretical support for these conditions. Second, we analyze dataset size requirements for AR-SA data-driven approximation based on the reduced dataset, taking into account solution reliability. Third, we compare the performance of the AR-SA approach with SA constructed on reduced scenario sets. The scenario reduction methods used are Fast Forward (FF) \cite{dupavcova2003scenario}, Simultaneous Backward (SB) \cite{heitsch2003scenario}, and $K$-Means \cite{keutchayan2023problem}. We use standard SA as a baseline. For the proposed AR-SA we observe nearly a twofold improvement in data efficiency compared to other scenario reduction techniques. We summarize the paper's workflow in Figure \ref{fig:workflow}.

The rest of this paper is organized as follows. Section~\ref{sec:setup} provides background and problem setup of the multi-stage high-voltage optimal power flow. Section~\ref{sec:chancecontrol} discusses the setup of the multi-step high-voltage OPF with automated generation control. Section~\ref{sec:chancecontrol} formulates the chance-constrained problem under consideration. Section~\ref{sec:apriori} presents the sketch of the a-priori sample reduction approach, formalizes it, and proves its validity. Furthermore, we prove that AR-SA (reduced data-driven approximation) theoretically requires fewer data samples to produce a reliable solution than classical SA. Section~\ref{sec:emp} compares AR-SA with classical Monte Carlo-based SA and other scenario reduction methods such as Fast Forward (FF), Simultaneous Backward (SB), and the clustering K-Means method.
Finally, the conclusion is in Section \ref{sec:conclusion}.
\section{Background and Problem Setup}\label{sec:setup}
\subsection{DC Optimal Power Flow}
The high-voltage DC model is a widely used load flow model in power systems. 
Let $G = (V, E)$ be a power system graph with the set of $n$ nodes (buses) $V$ and the set of $m$ lines (edges) $E$; 
$p \in \mathbb{R}^{n_g}$, $p_d \in \mathbb{R}^{n_d}$, and $\theta\in\mathbb{R}^{n}$ be vectors representing power generations, demands, and phase angles, respectively. 
The system is balanced so the sum of all power injections is zero, $\sum_{i \in V} p^i = \sum_{i \in V} p^i_d$.
For clarity, we designate one bus as the slack bus, with its phase angle set as $\theta_s = 0$. 
The components of admittance matrix $B$, $B \in \mathbb{R}^{n \times n}$, denoted as $B^{ij}$, are non-zero if there is a line between buses $i$ and $j$. For each node $i$, $B^{ii}$ is defined as the negative sum of the off-diagonal elements $B^{ij}$ with $j \neq i$. 
The DC power flow equations, security and constraints for $i \in V, (i,j)\in E$ are:
\vspace{-0mm}
\begin{gather*}
p-p_{d} = B \theta, \!\sum_{i=1}^{n_g} p^i = \!\!\!\sum_{i=1}^{n_d} p^i_d, 
\underline{p}^i_g \leq p^i_g \leq \overline{p}^i, |\theta^i - \theta^j| \leq \bar{\theta}^{ij}
\vspace{-2mm}
\end{gather*}
{\color{black}{}The DC Optimal Power Flow (DC-OPF) feasibility set is linear, defined by voltage phases and power generation within the grid topology, line characteristics (admittance matrix), and power demand \cite[Chapter 4.1.4]{wood2013power}. The objective is to find an economically optimal active power generation profile across available generators while adhering to technical limits and system demand, which establish rules and constraints on power transfer throughout the system.}
The feasibility set of DC-OPF can be reformulated as a polytope $P = {p: Wp \leq b}$ in the vector space of active power generations $p \in \mathbb{R}^{n_g}$, where $W \in \mathbb{R}^{J \times n_g}$ and $b \in \mathbb{R}^J$. Here, $n_g$ denotes the number of controllable generators, and $J$ represents the number of constraints \cite{lukashevich2021importance, lukashevich2021power}. Reliability constraints are violated when the power generation vector $p$ falls outside the polytope $P$.
Solving DC-OPF consists of finding the power flow by minimizing a convex cost of power generation $c(p_g)$ subject to the constraints defined by $P$. 
\vspace{-0mm}
\subsection{Source of uncertainty and AGC}
\label{sec:fluctuations}
\vspace{-0mm}
The fluctuations affect the power balance in the system and are typically managed through primary and secondary control \cite{machowski2020power}. In this paper, we consider linear Automatic Generation Control (AGC). The AGC recourse adjusts the generation to a new setpoint $p^{t+1} = p^t + \alpha \xi^t$ \cite{roald2017chance,baros2021examining,mezghani2020stochastic} with $p^t \in \mathbb{R}^{n_g}$, $\alpha \in \mathbb{R}^{n_g}$, and $\xi^t \in \mathbb{R}$ representing the \emph{total} demand-generation imbalance.
The participation factors $\alpha$ for secondary control can slightly vary, enabling long-term grid stability and fast control~\cite{machowski2020power}. 
The total system imbalance $\xi^{t}$ represents the sum of power fluctuations due to the intermittency of renewable generation, demand instability, and intra-day electricity trading. Specifically, in a power system with nodes $\mathcal{B}$, $\xi^t$ is calculated as $\xi^t = \sum_{b \in \mathcal{B}} (\xi_b^d)^t - (\xi_b^g)^t$, where $(\xi_b^d)^t$ and $(\xi_b^g)^t$ are random variables modeling demand and generation fluctuations at bus $b$ at time $t$, respectively. These variables follow various distributions based on the source. For example, $(\xi_b^g)^t$ follows a beta distribution at a bus with a PV generator \cite{wang2010probabilistic}, and wind speed near a bus with a wind farm may follow a Weibull distribution, though power output distribution may vary between turbines \cite{dhople2012framework}. Although $\xi^t$ as the sum of these independent random variables, often follows a Gaussian distribution due to the Lyapunov or Lindeberg-Feller Central Limit Theorem \cite{scholz2011central, rouaud2013probability, draper2021practical}. We validated the Gaussianity assumption using the Shapiro-Wilk normality test on time series data from the RTS-GLMC project, which includes demand, hydro, rooftop PV, and wind farm generations \cite{barrows2019ieee}. The test indicates that the hypothesis of normal distribution is rejected only for August and November at a significance level of $\alpha=0.05$, which supports the validity of the Gaussian assumption for $\xi^t$. We further consider the stacked temporal uncertainty vector $\xi = (\xi^1, \dots, \xi^T) \sim \mathcal{N}(\mu, \Sigma)$, where $\Sigma$ models temporal correlations.
{\color{black}{}We consider a system influenced by fluctuations and formulate an optimization problem to determine an economically optimal control strategy and initial system setpoint that satisfies JCC on technical limits and demand.}
Table~\ref{tab:notation} summarizes the paper's notation. Upper indices denote elements of vectors and matrices. When clear from context, we use $\mathbb{P}$, $\mathbb{E}$, and $\mathbb{V}$ for probability, expectation, and variance without specifying the distribution.

\begin{table}[t]
    \centering
    \caption{Paper notation.}
    \begin{tabular}{|p{0.35cm}|p{2.95cm}|p{0.95cm}|p{2.9cm}|}
        ${P}$ & DC-OPF feas. set & $P_{JCC}$ & JCC DC-OPF feas. set \\
        $p$ & vector of generation & $n_g$ & \# of controllable gen. \\
        $\alpha$ & participation factors & $J$ & \# of constraints in $P$ \\
        $T$ & \# of modeling timestamps & $\mathcal{N}(\mu, \Sigma)$ & Gaussian distribution \\  
        $\xi^t$  & power balance mismatch $t$ & $\Phi$ & Standard Gaussian CDF \\
        $I_n$ & identity matrix & $R$ & ramp-up/down limits \\
        $\mathcal{P}_r$ & theoretical redundancy set & $\hat{\mathcal{P}}_r$ & sufficient redundancy set\\
    \end{tabular}
    \label{tab:notation}    
\end{table}


\section{Chance constraint multi-stage control}
\label{sec:chancecontrol}
In this section, we introduce the JCC discrete-time dynamic DC-OPF with AGC and outline a data-driven approximation leading to a solution feasible for the original chance-constrained problem with high probability.

\subsection{Chance constrained optimization}
\vspace{-0mm}

Consider a dynamical system with $T$, $T<\infty$, timestamps, and $\xi^t$, $1 \le t\le T$ - total power mismatch due to uncertainties at timestamp $t$. 
Let individual uncertainties follow a Gaussian distribution: $\xi^t \sim \mathcal{N}(0, (\sigma^t)^2)$, $1\le t\le T$, so that $\xi \sim \mathcal{N}(0, \Sigma), ~\Sigma \in \mathbb{R}^{T\times T}$ with marginals distributed as $\mathcal{N}(0, (\sigma^t)^2)$. The temporal binding between system timestamps is modeled through the ramp rates of generators, ensuring realistic rates of change in power outputs as $|p^t_i - p^{t-1}_i| \leq R_i$, where $R_i > 0$. The discrete-time dynamic chance-constrained optimization problem is then: 
\vspace{-2mm}
\begin{align}
        & \hspace{32mm} \min_{p^t, \alpha} \mathbb{E} \sum_{t=1}^T c(p^t) \label{eq:optimal_control}, \qquad \texttt{s.t.:} 
        \\ 
        & \; \mathbb{P} 
        \begin{pmatrix}
                Wp^t \leq b, p^t = p^{t-1} + \alpha\xi^t, |p_k^t - p_k^{t-1}| \leq R_k,\\
                 1 \leq k \leq n_g, ~1\leq t \leq T
        \end{pmatrix} \geq 1 - \eta.\nonumber
\end{align}
where $\eta\in (0, 1/2]$, $\mathbb{P}$ is a joint measure induced by the uncertainty distribution and $\alpha \in \mathbb{R}^{n_g}$ is participation factors. 
A compact statement of the Problem \eqref{eq:optimal_control} is:
\vspace{-3mm}
\[\min_{p^t, \alpha} \mathbb{E} \sum_{t=1}^T c(p^t) \]
\vspace{-3mm}
\begin{equation}
    \begin{aligned}
        \!\!\texttt{s.t.:}  & \mathbb{P}\!\! 
        \begin{pmatrix}
                \mathcal{W}^p p^0 + E^\tau \mathcal{W}^{\alpha} \cdot \alpha \leq \beta, 0 \leq \tau \leq T 
        \end{pmatrix}\!\geq\!1 - \eta,\!\!\!
    \end{aligned}
    \label{eq:optimal_control_2} 
\end{equation}

where $E^\tau \in \mathbb{R}^{J+2n_g \times J + 2n_g}$ is a diagonal matrix with first $J$ diagonal elements equal $(1^\tau) ^\top \xi$, the rest are $(e^\tau)^\top\xi$. Here
$1^{\tau} \in \mathbb{R}^T$ has components $1^{\tau}_i = 1, ~ 0 \leq i \leq \tau, ~ 0$ otherwise. The vector $e^\tau_i = 1, ~ i=\tau, e^\tau_i = 0$ in the other case. Later in the paper, given specific $i:1 \leq i \leq J + 2n_g$ and $1\leq \tau \leq T$ we refer the second term components as $(E^\tau_i)^\top\xi \cdot (\omega_i^\alpha)^\top \alpha$, where $E^\tau_i$ is $1^\tau$ for $i \leq J$ and $e^\tau$ for $i > J$.
Matrices are obtained as vertical stacks: $\mathcal{W}^p = \left( W^\top, 0, 0 \right)^\top$ and $\mathcal{W}^\alpha = \left(W^\top, I_{n_g}, -I_{n_g} \right)^\top$, $\mathcal{W}^p, ~ \mathcal{W}^\alpha \in \mathbb{R}^{(J + 2 \cdot n_g)\times n_g}$. The right hand side of Eq.~\eqref{eq:optimal_control_2}, $\beta = \left(b^\top, R, R \right)^\top$ with $R = \{R_k\}_{k=1}^{n_g}$ being the vector of ramp up/down limits. 
We use subscript to refer to rows of the matrices $\mathcal{W}^p, \mathcal{W}^\alpha, E^\tau$; $\omega_i^{p}, ~ \omega_i^{\alpha}$ for the rows of matrices $\mathcal{W}^p, \mathcal{W}^\alpha$. 

\subsection{Scenario approximation of chance constrained control}
A Scenario Approximation (SA) of the Problem~\eqref{eq:optimal_control_2} via the set of scenarios $\xi(j), ~ j=1,\dots, N$, implying a separate set of constraints for each one, is: 
\vspace{-2mm}
    \begin{align}
        & \qquad \min_{p^0, \alpha} c(p^0), \qquad \texttt{s.t.:}  \label{eq:optimal_control_sampling_02} 
        \\ 
        \forall j, 1\leq j \leq N\!\!: & \;  \mathcal{W}^p p^0 + (E^\tau)^\top \xi(j) \mathcal{W}^{\alpha} \alpha \leq \beta, 0 \leq \tau \leq T.\nonumber
\end{align}
SA offers practical benefits but demands a huge number of samples to obtain a reliable solution \cite{calafiore2006scenario}. Reduction strategies, especially a-priori methods, are underexplored for JCC problems in power systems. Reducing data samples lowers constraint numbers, improving optimization tractability. We analyze sample redundancy conditions and provide reliability guarantees for the data-driven approximation's solution.
\section{A-priori scenario redundancy}
\label{sec:apriori}
{\color{black}{}
We define redundant data samples in an SA by setting analytical conditions for data redundancy in a JCC multi-timestamp DC-OPF. We provide the minimum number of scenarios needed for a $1 - \rho$ reliable solution, i.e., a solution feasible for Problem \eqref{eq:optimal_control_2} with a $1 - \rho$ probability, and derive reduction factors based on the measure of the analytical redundancy set.}
\subsection{Redundant scenarios}
{\color{black}{}In this subsection, we present an approach and theorems for a-priori scenario redundancy classification. First, we formalize data sample redundancy and illustrate it. Next, we provide formal statements that define an inner approximation of the redundancy set $\mathcal{P}_r$, classifying data samples as redundant.}
\begin{definition}
\label{def:redundant}
Let $\mathcal{I} = \{1, \dots, N\}$. Scenarios indexed with $\mathcal{I}_r \subset \mathcal{I}$ are called redundant iff a solution of SA with constraints corresponding to scenarios indexed with $\mathcal{I}_r$ are omitted - $(\hat{p}^0_{\mathcal{I} \setminus \mathcal{I}_r}, \hat{\alpha}_{\mathcal{I} \setminus \mathcal{I}_r})$ - is feasible for initial JCC \eqref{eq:optimal_control_2} and solution of SA $(\hat{p}^0_{\mathcal{I}_r}, \hat{\alpha}_{\mathcal{I}_r})$ with constraints corresponding to those scenarios indexed with $\mathcal{I}_r$ is not feasible for JCC \eqref{eq:optimal_control_2}.
\end{definition}
The redundancy concept is illustrated in Figure~\ref{fig:idea}. Our goal is to define conditions for a-priori identification of redundant scenarios. Consider a set $\mathcal{P}_r$ containing all redundant samples (scenarios indexed by $\mathcal{I}_{r}$ from Def.~\ref{def:redundant}). Keeping only these samples from $\mathcal{P}_r$ results in an infeasible solution for Problem \eqref{eq:optimal_control_2}. But keeping scenarios outside of $\mathcal{P}_r$ yields a feasible solution. Since it's challenging to analytically define $\mathcal{P}_r$, we aim to construct an inner approximation $\hat{\mathcal{P}}_r$. Such approximation provides a sufficient condition for identifying redundant samples. We derive the inner redundancy set $\hat{\mathcal{P}}_r$ for efficient scenario generation. Additionally, we employ a standard technical assumption \cite{campi2011sampling},  that ensures that problems with finitely many constraints are feasible:
\begin{assumption}\label{asmp:10}
For all possible uncertainty realizations $\xi(1), \dots, \xi(N)$, optimization problem \eqref{eq:optimal_control_sampling_02} is either infeasible or has a unique optimal solution.
\end{assumption}
\def\shift{2.7}
\def\shiftt{2.}
\def\shiftuselesssample{0.3}
\def\shiftusefulsample{1.1}
\def\xx2{0 + \shiftt}
\def\yx2{0 - \shiftt}
\def\threshold{33} 
\tikzset{snake it/.style={decorate, decoration={coil,amplitude=1pt, segment length=9pt}}}
\begin{figure}
\begin{minipage}{0.35\linewidth}
    \centering
        \begin{tikzpicture}[scale=0.45, node distance={15mm}, thick, main/.style = {draw, scale=.5}] 
        \clip (0,2) rectangle + (6,-7);
        
        \coordinate (x2) at (\xx2, \yx2);
        \coordinate (x2origin) at (\xx2 + 1.5, \yx2 - 2);
        \draw[olive] (x2) circle (1pt);
        \draw (x2) ++(1.1, -2.3) node[olive, above right] (tmp) {$(\hat{p}_{\mathcal{I}}, \hat{\alpha}_{\mathcal{I}})$};
        \draw[->, olive] (x2origin) -- (x2);
        \coordinate (a) at ( 4.755282581475767 , 1.545084971874737 );
        \coordinate (b) at ( 3.061616997868383e-16 , 5.0 );
        \coordinate (c) at ( -4.755282581475767 , 1.5450849718747375 );
        \coordinate (d) at ( -2.9389262614623664 , -4.045084971874736 );
        \coordinate (e) at ( 2.9389262614623646 , -4.045084971874738 );
        \draw (a) -- (b);
        \draw (b) -- (c);
        \draw (c) -- (d);
        \draw (d) -- (e);
        \draw (e) -- (a);
        \draw
        (a) ++(0.1, -1.21) node[below] (tmp) {$P$};
        
        \coordinate (ag) at ( 3.804226065180614 , 1.2360679774997896 );
        \coordinate (bg) at ( 2.4492935982947064e-16 , 4.0 );
        \coordinate (cg) at ( -3.804226065180614 , 1.23606797749979 );
        \coordinate (dg) at ( -2.351141009169893 , -3.2360679774997894 );
        \coordinate (eg) at ( 2.3511410091698917 , -3.2360679774997902 );
        \draw [black, snake it] (ag) -- (bg);
        \draw [black, snake it] (bg) -- (cg);
        \draw [black, snake it] (cg) -- (dg);
        \draw [black, snake it] (dg) -- (eg);
        \draw [black, snake it] (eg) -- (ag);
        \draw
        (ag) ++(0., -0.6) node[left, black] (tmp) {${P}_{\textup{JCC}}$};
        
        \coordinate (ar2) at ( 0.9510565162951535 + \shiftt , 0.3090169943749474 - \shiftt);
        \coordinate (br2) at ( 6.123233995736766e-17 + \shiftt , 1.0  - \shiftt);
        \coordinate (cr2) at ( -0.9510565162951535  + \shiftt, 0.3090169943749475  - \shiftt);
        \coordinate (dr2) at ( -0.5877852522924732 + \shiftt, -0.8090169943749473 - \shiftt);
        \coordinate (er2) at ( 0.5877852522924729 + \shiftt, -0.8090169943749476  - \shiftt);
        \draw [teal, snake it] (ar2) -- (br2);
        \draw [teal, snake it] (br2) -- (cr2);
        \draw [teal, snake it] (cr2) -- (dr2);
        \draw [teal, snake it] (dr2) -- (er2);
        \draw [teal, snake it] (er2) -- (ar2);
        \draw
        (br2) ++(1em, .1em) node[above, teal] (tmp) {$\mathcal{P}_{r}$};

        \coordinate (anvc) at ( 0.9510565162951535 * 0.7 + \shiftt , 0.3090169943749474 * 0.7 - \shiftt);
        \coordinate (bnvc) at ( 6.123233995736766e-17 + \shiftt , 1.0 * 0.7  - \shiftt);
        \coordinate (cnvc) at ( -0.9510565162951535 * 0.7  + \shiftt, 0.3090169943749475 * 0.7  - \shiftt);
        \coordinate (dnvc) at ( -0.5877852522924732 * 0.7 + \shiftt, -0.8090169943749473 * 0.7 - \shiftt);
        \coordinate (envc) at ( 0.5877852522924729 * 0.7 + \shiftt, -0.8090169943749476 * 0.7  - \shiftt);
        \draw [purple] (anvc) -- (bnvc);
        \draw [purple] (bnvc) -- (cnvc);
        \draw [purple] (cnvc) -- (dnvc);
        \draw [purple] (dnvc) -- (envc);
        \draw [purple] (envc) -- (anvc);
        \draw
        (envc) ++(3.3em, .3em) node[above, purple] (tmp) {$\hat{\mathcal{P}}_{r}$};

        
        \pgfmathsetmacro{\xRange}{1.3} 
        \pgfmathsetmacro{\yRange}{1.3} 
        \pgfmathsetseed{3}
        \foreach \i in {1,...,30} {
            
            \def\xrandom{\shiftt + rand*\xRange}
            \def\yrandom{-\shiftt + rand*\yRange}
        
            \filldraw [black] (\xrandom, \yrandom) circle (1pt);
            
        }
        \end{tikzpicture}  
        \end{minipage}
            \hfill 
        \begin{minipage}{0.65\linewidth}
        \vspace{-3mm}
            \begin{itemize}
                \item $P$ is a feasibility set
                \item $P_{\textup{JCC}}$ is a JCC feasibility set
                \item The black dots -- potential setpoints achievable by the AGC due to power fluctuations
                \item $(\hat{p}_{\mathcal{I}}, \hat{\alpha}_{\mathcal{I}})$ is the SA solution based on all data samples
            \end{itemize}
        \end{minipage}
        \vspace{-4mm}
    \caption{
    All data samples can be divided into redundant and non-redundant, depending on whether they are inside or outside the set $\mathcal{P}_r$ with an unknown structure. In practice, one can derive inner approximations $\hat{\mathcal{P}}_r$ of $\mathcal{P}_r$. The latter can be used to classify data by redundancy in SA.}
    \label{fig:idea}
\end{figure}

To get an analytical sufficient condition on the redundant scenarios, we derive a necessary condition for JCC feasibility converted to a sufficient condition of not being in the JCC feasibility set. Note that there is the following bound on feasibility probability from \eqref{eq:optimal_control_2}:
\begin{lemma}
    \label{lemma:bound_prob}
    {\color{black} Let} $\pi(p^0, \alpha) \geq 1 -\eta$, where $\pi(p^0, \alpha) = \mathbb{P}\left\{ \cap_{i, \tau} (\omega^p_i)^\top p^0 + (E^\tau_i)^\top \xi \cdot(\omega^\alpha_i)^\top \alpha \leq \beta_i\right\}$. Then 
    \vspace{-2mm}
    \[\max_{i, t} \mathbb{P} \left\{ (\omega^p_i)^\top p^0 + (E^\tau_i)^\top \xi \cdot (\omega^\alpha_i)^\top \alpha > \beta_i \right\} \leq 1-\pi(p^0, \alpha).\]
\end{lemma}
\begin{proof}
    As $\mathbb{P}\left\{ \cup_{i, t} (\omega^p_i)^\top p^0 + (E^\tau_i)^\top \xi \cdot (\omega^\alpha_i)^\top \alpha > \beta_i \right\} = 1 - \pi(p^0, \alpha)$, applying the Boole-Fréchet bound \cite[Theorem 4.2.1]{williamson1989probabilistic} to the left-hand side of this we prove the lemma.
\end{proof}
Lemma \ref{lemma:bound_prob} provides a handful necessary condition:
\begin{corollary}
    \label{lemma:corollary}
    Let $(p^0, \alpha)$ be feasible to JCC in \eqref{eq:optimal_control_2}. Then $\max_{i, t} \mathbb{P} \left\{ (\omega^p_i)^\top p^0 + (E^\tau_i)^\top \xi \cdot (\omega^\alpha_i)^\top \alpha > \beta_i \right\} \leq \eta$.
\end{corollary}
\begin{proof}
    Feasibility yields $\pi(p^0, \alpha) \geq 1 -\eta$ iff $1 - \pi(p^0, \alpha) \leq \eta$. Applying Lemma \ref{lemma:bound_prob} proves the corollary.
\end{proof}

We now formalize $\hat{\mathcal{P}}_r$ for Problem~\eqref{eq:optimal_control_2}. This gives us a-priori sufficient condition on sample redundancy.

\begin{theorem}
Let scenarios $\xi(j) \sim \mathcal{N}(0, \Sigma), ~ j \in \mathcal{I}=\{1, \dots, N\}$ form SA \eqref{eq:optimal_control_sampling_02} and a solution of this problem $(\hat{p}^0_{\mathcal{I}}, \hat{\alpha}_{\mathcal{I}})$ be feasible for the JCC problem \eqref{eq:optimal_control_2}. Moreover, assume that the cost function $c(\cdot)$ is linear. {\color{black} Let $\hat{\mathcal{P}}_r = \{ \xi\in \mathbb{R}^T:~ |(E^{\tau}_i)^\top \xi|  \leq \Phi^{-1}(1 - \eta) \sigma^\tau_i \gamma ~\forall i, \tau \}$,} where $\gamma \in (0, 1)$ and $(\sigma^\tau_i)^2 = (E^\tau_i)^\top \Sigma (E^\tau_i)$.
Then, first, SA {\color{black} where scenarios} $\xi(j), ~ j \in \mathcal{I}_r = \{ j: \xi(j) \in \hat{\mathcal{P}}_r \}$ yields solution $(\hat{p}^0_{\mathcal{I}_r}, \hat{\alpha}_{\mathcal{I}_r})$ that is not feasible for original JCC Problem \eqref{eq:optimal_control_2}. Second, SA {\color{black} where scenarios} $\xi(j), ~ j \in \mathcal{I} \setminus \mathcal{I}_r$ yields the solution $(\hat{p}^0_{\mathcal{I} \setminus \mathcal{I}_r}, \hat{\alpha}_{\mathcal{I} \setminus \mathcal{I}_r})$ that is feasible for the original JCC Problem \eqref{eq:optimal_control_2}.
\label{th:P_r sampling polytope}
\end{theorem}
\begin{proof}
    The feasibility set of SA is given by $(w^p_i)^\top p^0 + (E_i^\tau)^\top \xi(j) \cdot (\omega^\alpha_i)^\top \alpha \leq \beta_i, ~ \forall i, \tau, j$. Since the cost function is linear for solution of SA, $\exists i', \tau', j': ~ (w^p_{i'})^\top p_{\mathcal{I}_r}^0 + (E_{i'}^{\tau'})^\top \xi(j') \cdot (\omega^\alpha_{i'})^\top \alpha_{\mathcal{I}_r} = \beta_{i'}$. Next, one has {\color{black} $(E^{\tau}_i)^\top \xi (j')  \leq \Phi^{-1}(1 - \eta) \sigma^\tau_i \gamma$, because $
    \xi(j') \in \hat{\mathcal{P}}_r$, positive absolute value case}. This implies that $(\omega^p_{i'})^\top p^0_{\mathcal{I}_r} \geq \beta_{i'} - \Phi^{-1}(1-\eta) \sigma^{\tau'}_{i'} \| (\omega_{i'}^\alpha)^\top \alpha_{\mathcal{I}_r} \| \gamma$. However, the necessary condition from Corollary \ref{lemma:corollary} implies that $\forall i, \tau ~ (\omega_i^p)^\top p^0_{\mathcal{I}_r} \leq \beta_i - \Phi^{-1}(1-\eta) \sigma^\tau_i \| (\omega_i^\alpha)^\top\alpha_{\mathcal{I}_r}\|$. Recalling that $\gamma \in (0, 1)$ we obtain contradiction that leads to the fact that $(\hat{p}_{\mathcal{I}_r}^0, \hat{\alpha}_{\mathcal{I}_r})$ is infeasible for original JCC Problem.
    Now drop redundant scenarios from the SA. Since $(\hat{p}^0_{\mathcal{I}}, \hat{\alpha}_{\mathcal{I}})$ is feasible for JCC problem and data samples $\xi(j), ~ j \in \mathcal{I}_r$ do not contribute to the feasibility, then SA solution $(\hat{p}_{\mathcal{I} \setminus \mathcal{I}_r}^0, \hat{\alpha}_{\mathcal{I} \setminus \mathcal{I}_r})$ built on $\xi(j), ~ j \in \mathcal{I} \setminus \mathcal{I}_r$ is feasible for JCC problem.
\end{proof}

Theorem \ref{th:P_r sampling polytope} establishes a sufficient, i.e., a priori condition on sample redundancy across a given dataset $\xi(j), j \in \mathcal{I}$ that guarantees a feasible solution. Specifically, if a dataset provides assurance of producing a feasible solution for the original JCC problem, samples within $\hat{\mathcal{P}}_r$ may be disregarded. Thus this theorem offers a way to assess the dataset's potential a-priori: if all data samples fall within $\mathcal{P}_r$, it is impossible to derive any feasible solution for the JCC from this data.

\subsection{SA Solution Guarantees and Dataset Complexity}
In this section, we provide theoretical guarantees on the solution of reduced SA that does not contain redundant samples within $\hat{\mathcal{P}}_r$.
The following theorem addresses the sampling complexity of the SA with data samples indexed $\mathcal{I} \setminus \mathcal{I}_r$.
\begin{theorem}\label{thm:40}
Let $(\hat{p}^0, \hat{\alpha})$ be a unique solution of the SA Problem~\eqref{eq:optimal_control_sampling_02} with $N$ i.i.d. samples, so that none of the samples belong to $\hat{\mathcal{P}}_{r}$. Moreover, assume that for any $N$ the assumption \ref{asmp:10} is fulfilled. Then for any $\rho \in (0,1)$ and any~$\eta \in (0, 1/2]$, $(\hat{p}^0, \hat{\alpha})$ is also a solution for the Chance-constrained optimal power flow Problem~\eqref{eq:optimal_control_2} with probability at least $1-\rho$ if 
$
  N \ge \left\lceil 2\eta^{-1}(1-\nu)\ln \frac{1}{\rho} + 2d + 2d\eta^{-1}(1-\nu) \ln\frac{2(1-\nu)}{\eta} \right\rceil, 
$
 where $d$ is a dimension of the space of controllable generators and participation factors, i.e., $d = 2 n_g$, and $\nu$ is the probability of a random scenario $\xi \sim \cN(0, \Sigma)$ to belong to $\hat{\mathcal{P}}_{r}$, and $\nu < 1$. 
\end{theorem}
\begin{proof}
First, notice that discarding random scenarios $\xi \notin \hat{\mathcal{P}}_{r}$ is equivalent in solving the SA problem with sampling $\xi$ from a distribution $\mathcal{D}$ where  $\xi \sim \mathcal{D} \Leftrightarrow \xi\sim \cN(0, \Sigma) \st \xi\not\in \hat{\cP}_{r}$. From the theorem statement, $1-\nu$ is the probability mass associated with samples in $\xi\sim \cN(0, \Sigma)$ that are outside $\hat{\cP}_{r}$.
{\color{black} Assumption \ref{asmp:10} and convexity of each function in Problem \ref{eq:optimal_control_2} meet the conditions of Calafiori and Campi~\cite[{\color{black}Theorem 1}]{calafiore2006scenario} implying that} for any probability $\rho \in (0,1)$ and any confidence threshold probability $\varepsilon$, and dimension of the space of parameters $d$ one has, for~$N_1$
\vspace{-1mm}
\begin{align}\label{eq:mon}
\vspace{-2mm}
  N_1 \ge \left\lceil 2\varepsilon^{-1} \ln (1/\rho) + 2d + 2d\varepsilon^{-1} \ln (2/\varepsilon)\right\rceil 
\end{align}
scenarios from $D$ and the optimal solution 
$(\hat{p}^0, \hat{\alpha})$ of the Problem~\eqref{eq:optimal_control_sampling_02}, 
the probability of failure is bounded as $
  \mathbb{P}_D\{(\hat{p}^0, \hat{\alpha}) \textup{ is feasible for \ref{eq:optimal_control_2}}\} \geq 1-\varepsilon
$ with prob. at least $1-\rho$. 

Notice, that the bounds on the number of samples (see Eq.~\eqref{eq:mon}) is strictly decreasing in $\varepsilon$ for $\varepsilon \in (0,1)$. As scenarios in $\hat{\cP}_{r}$ are redundant and do not contribute to overall solution reliability, to get a probability of failure $\eta$ according to measure $\cN(0, \Sigma)$, we need the failure probability according to $\mathcal{D}$ to be at least $\varepsilon = \eta/(1-\nu)$.  
Thus, using $\varepsilon = \eta/(1-\nu)$ and monotonicity of~Ineq.~\eqref{eq:mon} one completes the proof. 
\end{proof}

\section{Empirical Study}\label{sec:emp}
\subsection{Algorithms and Implementation Details}
We compare the performance of SA based on different strategies: classical Monte-Carlo (SA) and the proposed A-priori Reduced (AR-SA), and 
the state-of-the-art {\color{black} scenario reduction methods: Fast-Forward, Simultaneous Backward \cite{heitsch2003scenario, dupavcova2003scenario} and K-Means \cite{keutchayan2023problem}}. 
For our test cases, we consider power systems from MATPOWER \cite{zimmerman2010matpower}, specifically Grid-6WW \cite{wood2013power} (pp. 104, 112, 119, 123-124, 549), Washington-14 and IEEE-30.
We implemented the algorithms using Python 3.9.13 and PandaPower 2.8.0 \cite{pandapower.2018} on a MacBook Pro (M1 MAX, 64 GB RAM). The optimization problems were solved using the Pyomo framework \cite{hart2017pyomo}
, which employed GLPK \cite{Oki2012GLPKL}.
The code is available on GitHub\footnote{https://github.com/vjugor1/OptimalControlScenarioApproximation}.
\subsection{Test Cases and Numerical Results}
We conducted two case studies to evaluate the performance of SA, AR-SA, and other scenario reduction methods under different scenarios. {\color{black} The first study focused on the Grid-6WW and Washington-14, comparing the number of samples $N$ needed to achieve the $1-\rho=0.99$ reliability of $1-\eta$ feasible solution among 5 methods}. 
{\color{black}The second study analyzed IEEE-30 bus system, which consist of 30 buses. 
In this case, we compared the number of samples needed to achieve solution reliability of $1-\rho=0.99$ and assessed the total execution time, including the scenario reduction step. We summarized the required number of samples in Table \ref{tab:summary_results}.}
In all case studies, we model power generation and consumption fluctuations with a standard deviation of 0.01 of their nominal values, increasing cumulatively for each temporal snapshot. Thus, we expressed fluctuations as $\xi^t \sim \mathcal{N}(0, (\tilde{\sigma_t}^2))$, where $\tilde{\sigma_t}^2 = \tilde{\sigma}_0^2 \cdot t$. Simulations were carried out over $T=3$ temporal snapshots.


{\bf Evaluation methodology.}
To get an empirical estimation of the solution reliability $1-\hat{\rho}$, we independently construct $L=100$ different approximations for both SA and reduced SA. Further, we short name SA with scenarios reduced by FF, SB, and K-Means and $\hat{\mathcal{P}}_r$ as SA-FF, SA-SB, SA-KMeans, and AR-SA respectively. The data-driven approximation size $N$ starts with $3$ and increases until the corresponding scenario approximation problem reaches $1-\hat{\rho} > 0.99$. For each approximation, we obtain $L$ different solutions: $(x^*_N)_l, ~ l=1, \dots, L$. We estimate the confidence of each obtained solution by running out-of-sample validation. We use $10^4$ Monte-Carlo samples of $\xi$ to estimate of JCC feasibility constraint $(\hat{\mathbb{P}}_N)_l$, $l=1, \dots, L$. Finally, the solution reliability is given by $1 - \hat{\rho}$, which represents the fraction of $L$ solutions $(x^*_N)_l$ such that $(\hat{\mathbb{P}}_N)_l \geq 1 - \eta$. 
Alg.~\ref{alg:estimate_delta} summarizes the sequence of steps.

{\normalsize
\begin{algorithm}[!t]
\caption{$\hat{\rho}$ -- an empirical estimate}\label{alg:estimate_delta}
\begin{algorithmic}
\Require $L$ - $\#$ trials, DC-OPF parameters, $\eta$ - confidence level, $N_0$ - initial size of SA, $N_{\max}$ - maximal size of SA
\State $N \gets N_0$, \; $\hat{ \boldsymbol \rho}$ -- storage for $\hat{\rho}_N$
\While{$N \leq N_{\max}$}
    \State $C_N \gets 0$ -- feasibility counter; \;$l \gets 1$
    \While{$l \leq L$}
        \State Obtain SA solution \eqref{eq:optimal_control_sampling_02} with $N$ samples
        \State Estimate feasibility prob. $(\hat{\mathbb{P}}_N)_l$ using MC. \label{alg:estimate_delta:phat_N_l}
        \If{$(\hat{\mathbb{P}}_N)_l \geq 1 - \eta$}
            $C_N \gets C_N +1$
        \EndIf
    \EndWhile
    \State $1-\hat{\rho}_N \gets C_N / L$ -- fraction of trials that are feasible
    \State Append $\hat{\rho}_N$ to $\hat{ \boldsymbol \rho}$;\; update  $n  \gets n + N_{\max}/ 10$
\EndWhile
\State \Return mean of $\hat{ \boldsymbol \rho}$
\end{algorithmic}
\end{algorithm}
}

\begin{figure*}[hbt]
\vspace{-0mm}
\begin{minipage}[b]{.7\textwidth}
\begin{subfigure}{.47\textwidth}
  \centering
  \hspace{-1mm}\includegraphics[width=0.9\linewidth]{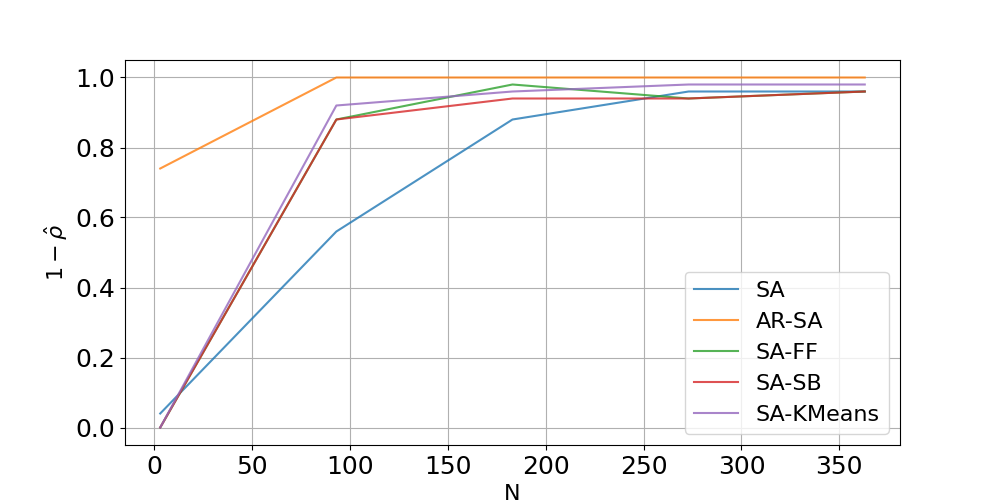}
  \caption{Empirical reliability 
  vs. $\#$ samples in CC-OPF for Washington 14 bus, $\eta = 0.01$.}
  \label{fig:washington14conservatism}
\end{subfigure}
\begin{subfigure}{.47\textwidth}
  \centering
  \hspace{-1mm}\includegraphics[width=0.9\linewidth]{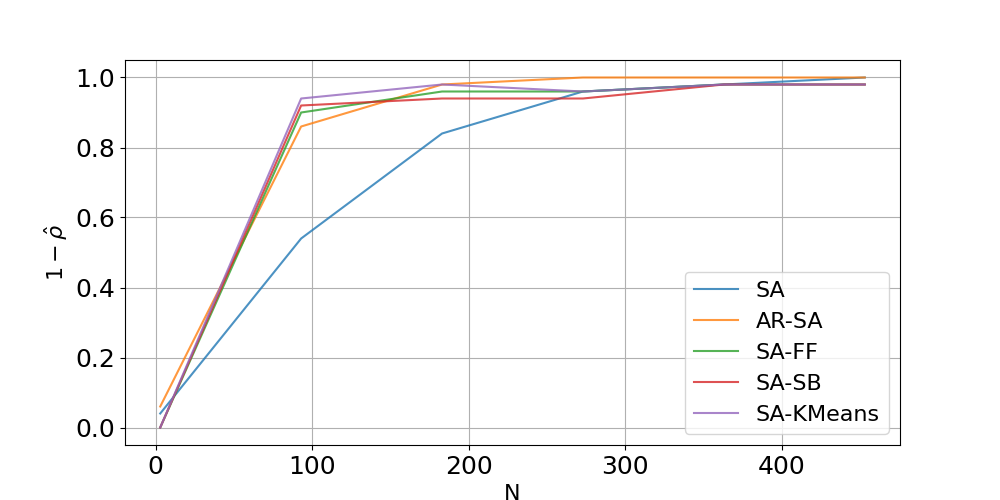}
  \caption{Empirical reliability 
  vs $\#$ samples  in CC-OPF for IEEE 30 bus system, $\eta = 0.01$.}
  \label{fig:ieee30conservatism}
\end{subfigure}\\
\begin{subfigure}{.47\textwidth}
  \centering
  \hspace{0mm}\includegraphics[width=0.9\linewidth]{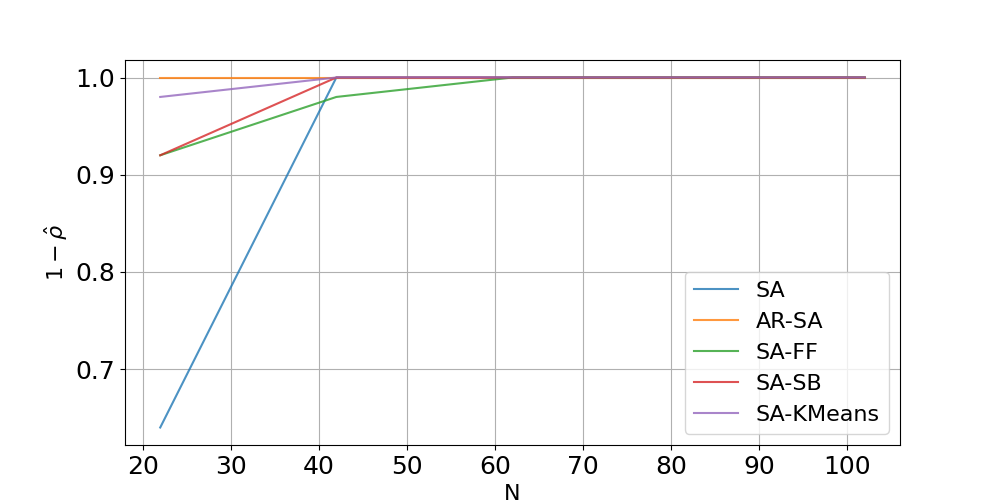}
  \caption{Empirical reliability 
  vs $\#$ samples  in CC-OPF for Grid6-WW 6 bus system, $\eta = 0.1$.}
  \label{fig:grid6reliability}
\end{subfigure}
\begin{subfigure}{.47\textwidth}
  \centering
\hspace{10mm}\includegraphics[width=0.9\linewidth]{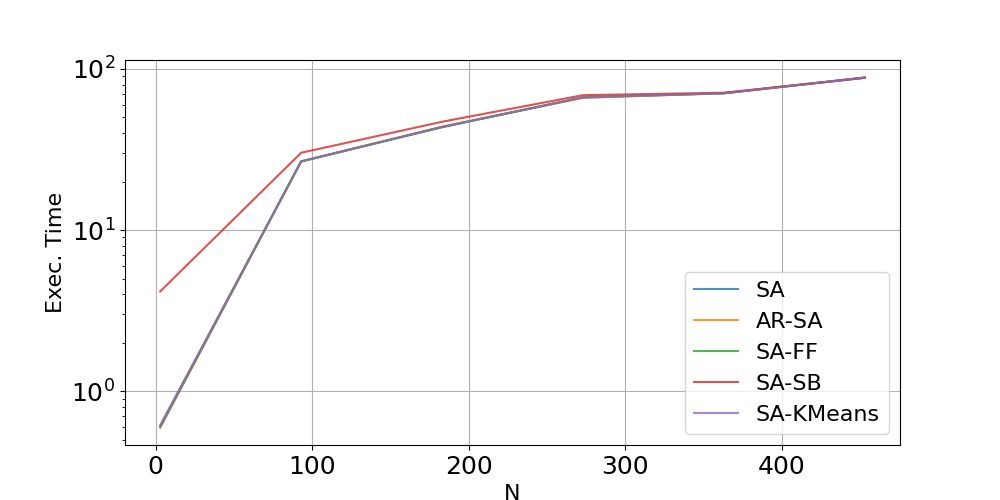}~~~~~~\hfill
  \caption{Exec. time vs $\#$ samples, IEEE 30 bus, $\eta = 0.01$.
  }
  \label{fig:ieee30time}
\end{subfigure}
\end{minipage}
\vspace{-1mm}
\begin{minipage}[b]{.3\textwidth}
\vspace{-1mm}
\caption{(a, b) Empirical reliability ($1-\hat{\rho}$). (c) Execution time for Grid6-WW. (d, e, f) the probability  of constraints feasibility (25\% -- 75\% interquartile range, IQR) $(\hat{\mathbb{P}}_N)_l$ for CC-OPF ($1-\eta =.99$) in the Washington 14 \& IEEE-30, Grid-6WW   systems. Circles stand for samples outside of the $\pm 1.5*\text{IQR}$. The empirical estimates are computed with $L = 200$ optimization instances (for $1-\hat{\rho}$), and $N_{MC}=10^4$ Monte-Carlo samples for each instance to determine constraint validation (for box-plot of $(\hat{\mathbb{P}}_N)_l$), as described in Algorithm~\ref{alg:estimate_delta}. 
}
\end{minipage}
\label{fig:ieee118}
\vspace{-10mm}
\end{figure*}

{\bf Complexity and Execution Time.}
In addition to solving the SA problem, Eq.~\eqref{eq:optimal_control_sampling_02}, selecting scenarios requires computational effort. Standard Monte Carlo-based SA (denoted as SA) does not require scenario reduction, unlike other methods that need additional computations for scenario selection.
The Fast Forward method adds scenarios one by one based on probabilistic metrics (2-Wasserstein distance) and redistributes probabilities after each addition. Simultaneous Backward, on the other hand, removes scenarios using the same process. Denoting $N_r$ as the target number of scenarios after reduction, the complexities of these methods are $O(N_r^3 + N_r N^2)$ and $O(N_r^3 + N^3)$, respectively \cite{heitsch2003scenario, rujeerapaiboon2022scenario}. The K-Means algorithm is based on iterative estimation of scenario cluster centers and requires estimation of $L_2$ distance between scenarios and cluster centers. This algorithm requires $O(N_rN^2)$ \cite{pakhira2014linear}. Contrary, in the proposed AR-SA the reduction step is conducted via checking if a current sample is within $\hat{\mathcal{P}}_r$, thus, the complexity is $O(N)$. The construction of $\hat{\mathcal{P}}_r$ itself grows linearly with the number of deterministic constraints under probability measure in JCC of \eqref{eq:optimal_control_2}.
We analyze total execution time, including reduction and solving the corresponding reduced scenario approximation. The test case is IEEE-30 grid with a target JCC feasibility level of $1-\eta=0.99$. The results, shown in Figure \ref{fig:ieee30time}, indicate that the execution time is almost similar for all reduction methods with classical Monte Carlo SA, except SA-SB. This suggests that scenario reduction is computationally inexpensive compared to solving the optimization problem. Additionally, this experiment supports the proposed method's practicality.

\vspace{-0mm}
{\bf SA Solution Reliability.}
\begin{table}[t]
    \centering
    \adjustbox{width=\linewidth}{
        \begin{tabular}{lr|rlrll}
        \toprule
        Case & $\eta$ & SA & AR-SA & SA-FF & SA-SB & SA-KMeans \\
\midrule
grid14 & 0.05 & 93 & 48 & 48 & 138 & 48 \\
grid30 & 0.05 & 138 & 93 & 138 & 138 & 93 \\
grid14 & 0.01 & 363 & 93 & 363 & 363 & 363 \\
grid30 & 0.01 & 453 & 273 & 453 & 453 & 453 \\
\bottomrule
        \end{tabular}
    }
    \caption{The number of samples for AR-SA and SA required in CC-OPF with a confidence threshold of $1-\eta$ to get the empirical reliability of $1-\hat{\rho} = 0.99$. The value of $1-\eta$ is given by out-of-sample Monte Carlo; the empirical reliability is given by averaging over $L=100$ independent CC-OPF problem instances, as described in Algorithm~\ref{alg:estimate_delta}. 
    }
    \label{tab:summary_results}
\end{table}
Following this analysis, we evaluate methods on larger grids with high reliability requirements. 
This experiment seeks to find the number of samples sufficient for a solution of data-driven approximation to be feasible for original JCC with high probability for different $\eta$.
We estimate the number of samples required to reach confidence thresholds of $1-\eta = 0.95$ and $0.99$ with reliability of $0.99$ for AR-SA, SA, SA-FF, SA-SB and SA-KMeans in Table~\ref{tab:summary_results}, showing the number of samples required is 30-50\% less for AR-SA compared to classical SA and the advantage of AR-SA increases with the increase of $1-\eta$.
We illustrate the dependence between empirical reliability, $1-\hat{\rho}$, and the number of samples, $N$, for different values across the IEEE-30, Washington-14 and Grid6-WW systems in Figures \ref{fig:ieee30conservatism},  \ref{fig:washington14conservatism} and \ref{fig:grid6reliability}, respectively. We maintain a confidence threshold for JCC feasibility of $1-\eta = 0.99$. 
Notably, AR-SA achieves higher reliability levels ($1-\hat{\rho}$) with significantly fewer samples $N$. From Figures \ref{fig:ieee30conservatism},  \ref{fig:washington14conservatism} and \ref{fig:grid6reliability} one can observe that for Washington-14 and Grid6-WW AR-SA reaches high reliability levels with the lower number of samples and in IEEE-30 case, though SA, SA-FF, SA-SB, and SA-KMeans reach $1 - \rho = 0.9$ with less samples, the higher $1 - \rho =0.99$ reliability level is reached by AR-SA first. The latter is due to the problem-specific redundancy set $\hat{\mathcal{P}}_r$ construction that is able to filter specifically those scenarios that are irrelevant to the given problem.

\vspace{-1mm}
\section{Conclusion}\label{sec:conclusion}
\vspace{-1mm}
Data-driven approximations are useful in chance-constrained stochastic programs with unknown uncertainty distribution and/or JCC settings. However, the data requirements rapidly become infeasible with the increase of size and reliability requirements. To address this, we proposed a novel approach that allows to a-priorily identify and remove redundant scenarios in stochastic approximations for JCC dynamic multi-timestamp DC-OPF. We prove the validity of this approach theoretically and 
ensured its high empirical performance over various test cases. 
\vspace{-2mm}
\bibliographystyle{IEEEtran}
\bibliography{biblio}
\end{document}